\RequirePackage[l2tabu, orthodox]{nag}

\documentclass[12pt]{amsart}
\usepackage{fullpage,url,amssymb,enumerate,colonequals}
\usepackage{mathrsfs} 

\usepackage[OT2,T1]{fontenc}
\DeclareSymbolFont{cyrletters}{OT2}{wncyr}{m}{n}
\DeclareMathSymbol{\Sha}{\mathalpha}{cyrletters}{"58}

\usepackage{color}

\newcommand{\defi}[1]{\textcolor{blue}{\textsf{#1}}} 

\newcommand{\Q}{\mathbb{Q}}
\newcommand{\R}{\mathbb{R}}
\newcommand{\Z}{\mathbb{Z}}

\newtheorem*{theorem}{Theorem}

\theoremstyle{definition}
\newtheorem*{definition}{Definition}

\theoremstyle{remark}

\usepackage{microtype}

\usepackage[
	backref,
	pdfauthor={Bjorn Poonen}, 
]{hyperref}
\usepackage[alphabetic,backrefs,lite]{amsrefs} 

\begin{document}

\title{Why all rings should have a~1}
\author{Bjorn Poonen}
\thanks{The writing of this article was supported by National Science Foundation grant DMS-1069236.}
\address{Department of Mathematics, Massachusetts Institute of Technology, Cambridge, MA 02139-4307, USA}
\urladdr{\url{http://math.mit.edu/~poonen/}}
\date{April 1, 2014}


\maketitle

\section{Introduction}\label{S:introduction}

Should the definition of ring require
the existence of a multiplicative identity~$1$?

Emmy Noether, when giving the 
modern axiomatic definition of a commutative ring, in 1921, 
did not include such an axiom.\footnote{See \cite{Noether1921}*{p.~29}.  Noether was preceded by David Hilbert in 1897 and Adolf Fraenkel in 1914, who used the word \emph{ring} in more restrictive senses.  Hilbert used \emph{Zahlring}, \emph{Ring}, and \emph{Integr\"atsbereich} to mean what we would call a finitely generated subring of an algebraic closure of $\Q$ \cite{Hilbert1897}*{\S31}; he implicitly included a~$1$.  Fraenkel explicitly required a~$1$, but the concept he was axiomatizing was quite different from the modern concept of ring: for example, his axiom $R_{8)}$ required that every non-zerodivisor have a multiplicative inverse~\cite{Fraenkel1914}*{p.~144}.}
For several decades, algebra books followed suit.\footnote{See, for instance, \cite{VanDerWaerden1966}*{\S3.1} and \cite{Zariski-Samuel1975}*{I.\S5}.}
But starting around 1960, many books by notable
researchers\footnote{Examples include \cite{EGA-I}*{0.(1.0.1)}, \cite{Lang1965}*{II.\S1}, \cite{Weil1967}*{p.~XIV}, and \cite{Atiyah-Macdonald1969}*{p.~1}.}
began using the term ``ring'' to mean ``ring with~$1$''.
Sometimes a change of heart occurred in a single person, or between editions of a single book, always towards requiring a~$1$.\footnote{Compare \cite{Jacobson1951}*{p.~49} with \cite{Jacobson1985}*{p.~86}, or \cite{Birkhoff-MacLane1953}*{p.~370} with \cite{Birkhoff-MacLane1965}*{p.~346}, or \cite{BourbakiAlgebre1in1958}*{I.\S8.1} with \cite{BourbakiAlgebre1-3}*{I.\S8.1}.} 
Reasons were not given; perhaps it was just becoming 
increasingly clear that the~$1$ was needed for many theorems to hold.\footnote{Some good reasons for requiring a $1$ are explained in~\cite{Conrad-standard-defs}.}

But is either convention more \emph{natural}?
The purpose of this article is to answer yes,
and to give a reason: 
existence of a~$1$ is a part of what associativity should be.

\section{Total associativity}\label{S:associativity}

The whole point of associativity is that it lets us assign an
unambiguous value to the product of any finite sequence of 
two or more terms.
By why settle for ``two or more''?
Cognoscenti do not require \emph{sets} to have two or more elements.
So why restrict attention to sequences with two or more terms?
Most natural would be to require \emph{every} finite
sequence to have a product, even if the sequence is of length $1$ or $0$:

\begin{definition}
A \defi{product} on a set $A$ is a rule that assigns to each finite sequence
of elements of $A$ an element of $A$,
such that the product of a $1$-term sequence is the term.
A product is \defi{totally associative} 
if each finite product of finite products
equals the product of the concatenated sequence
(for example, $(abc)d(ef)$ should equal the 6-term product $abcdef$).
\end{definition}

The ring axioms should be designed so that they give rise
to a totally associative product.
Now the key point is the following theorem, 
the more involved direction of which (the ``if'' direction) is 
\cite{BourbakiAlgebre1-3}*{I.\S1.2, Th\'eor\`eme~1, and~\S2.1}.

\begin{theorem}
A binary operation extends to a totally associative product
if and only if it is associative\footnote{I.e., associative in the usual sense, on triples: $(ab)c=a(bc)$ for all $a,b,c \in A$.} 
and admits an identity element.
\end{theorem}

What?!  Where did that identity element come from?
The definition of totally associative implies the equations
\begin{align*}
	(abc)d &= abcd \\
	(ab)c &= abc \\
	(a)b &= ab \\
	()a &= a.
\end{align*}
The last equation, which holds for any $a$,
shows that the empty product $()$ is a left identity.
Similarly, $()$ is a right identity, so $()$ is an identity element.

Therefore it is natural to require rings to have a~$1$.
But occasionally one does encounter structures
that satisfy all the axioms of a ring except for the existence of a~$1$.
What should they be called?
Happily, there is an apt answer, suggested by Louis Rowen: 
\defi{rng}!\footnote{See \cite{Jacobson1985}*{p.~155}.  Bourbaki introduced a different pejorative for the same concept: \defi{pseudo-ring} \cite{BourbakiAlgebre1-3}*{I.\S8.1}.  Keith Conrad observes in \cite{Conrad-standard-defs}*{Appendix~A} that the usual definition of \defi{(associative) $\Z$-algebra} gives the same notion and that hence this terminology could be used.}
As our reasoning explains and as Rowen's terminology suggests, 
it is better to think of a rng as a ring with something missing
than to think of a ring with~$1$ as having something extra.

\section{Counterarguments}\label{S:rngs}

Here we mention some arguments for \emph{not} requiring a~$1$,
in order to rebut them.

\begin{itemize}
\item ``Algebras should be rings, but Lie algebras usually do not have a $1$.''

Lie algebras are usually not associative either.
We require a~$1$ only in the presence of associativity.
It is accepted nearly universally that ring multiplication should
be associative,
so when the word ``algebra'' is used in a sense broad enough to include
Lie algebras, it is understood that algebras have no reason to be rings.

\item ``An infinite direct sum of nonzero rings does not have a~$1$.''

Direct sums are typically defined for objects like vector spaces
and abelian groups, for which the set of homomorphisms between
two given objects is an abelian group, for which cokernels exist,
and so on.
Rings fail to have these properties, whether or not a~$1$ is required.
So it is strange even to speak of a direct sum of rings.
{}Category theory explains that the natural notion 
for rings is the direct \emph{product}.
Each ring may be viewed as an \emph{ideal} in the direct product;
then their direct sum is an ideal too.

\item ``If a~$1$ is required, 
then function spaces like the space $C_c(\R)$
of compactly supported continuous functions $f \colon \R \to \R$
will be disqualified.''

This is perhaps the hardest to rebut, given the importance of function spaces.
But many such spaces are ideals in a natural ring
(e.g., $C_c(\R)$ is an ideal in the ring $C(\R)$ of \emph{all} 
continuous functions $f \colon \R \to \R$),
and fail to include the $1$ only because of some condition imposed
on their elements. 
So one can say that they, 
like the direct sums above and like the rng of even integers, 
\emph{deserve} to be ousted from the fellowship of the ring.
In any case, however, these function spaces still qualify as $\R$-algebras.
\end{itemize}

\section{Final comments}\label{S:final comments}

Once the role of the empty product is acknowledged,
other definitions that seemed arbitrary become natural.
A ring homomorphism $A \to B$ should respect finite products,
so in particular it should map 
the empty product $1_A$ to the empty product $1_B$.
A subring should be closed under finite products,
so it should contain the empty product $1$.
An ideal is prime if and only if its complement 
is closed under finite products; in this case, $1$ is in the complement;
this explains why the unit ideal $(1)$ in a ring is never considered to 
be prime.

The reasoning in Section~\ref{S:associativity} involved 
only one binary operation,
so it explains also why monoids are more natural than semigroups.\footnote{A \defi{semigroup} is a set with an associative binary operation, and a \defi{monoid} is a semigroup with a~$1$.}
Similar reasoning explains why 
the axioms for a category require not only 
compositions of two morphisms but also identity morphisms: 
given objects $A_0,\ldots,A_n$ and a chain of morphisms
\[
	A_0 \stackrel{f_1}\to A_1 \stackrel{f_2}\to 
	\cdots \stackrel{f_n}\to A_n,
\]
one wants to be able to form the composition, even if $n=0$.

It would be ridiculous to introduce the definition of ring
to beginners in terms of totally associative products.
But it is nice to understand why certain definitions should be
favored over others.


\section*{Acknowledgement} 

I thank Keith Conrad for many helpful comments.

\end{document}